\documentclass[12pt]{article}
\usepackage{amsfonts}
\usepackage{hyperref}
\usepackage{mathrsfs}
\usepackage{indentfirst}
\usepackage{amsmath}
\usepackage{amssymb}
\usepackage[amsmath,thmmarks]{ntheorem}
\usepackage{cite}
\usepackage{graphicx}
\usepackage{tikz}

\newtheorem{definition}{\bf Definition}[section]
\newtheorem{lemma}{\bf Lemma}[section]
\newtheorem{theorem}{\bf Theorem}[section]

\newtheorem{example}{\bf Example}[section]
\newtheorem{algorithm}{\bf Algorithm}[section]

\begin{document}
\setcounter{page}{1}

\title{{\textbf{Eigenproblems in addition-min algebra}}\thanks {Supported by
the National Natural Science Foundation of China (No.12071325)}}
\author{Meng Li\footnote{\emph{E-mail address}: 1228205272@qq.com}, Xue-ping Wang\footnote{Corresponding author. xpwang1@hotmail.com; fax: +86-28-84761393},\\
\emph{School of Mathematical Sciences, Sichuan Normal University,}\\
\emph{Chengdu 610066, Sichuan, People's Republic of China}}

\newcommand{\pp}[2]{\frac{\partial #1}{\partial #2}}
\date{}
\maketitle
\begin{quote}
{\bf Abstract} In order to guarantee the downloading quality requirements of users and improve the stability of data transmission in a BitTorrent-like peer-to-peer file sharing system, this article deals with eigenproblems of addition-min algebras. First, it provides a sufficient and necessary condition for a vector being an eigenvector of a given matrix, and then presents an algorithm for finding all eigenvalues and eigenvectors of a given matrix. It further proposes a sufficient and necessary condition for a vector being a constrained eigenvector of a given matrix and supplies an algorithm for computing all the constrained eigenvectors and eigenvalues of a given matrix. This article finally discusses the supereigenproblem of a given matrix and presents an algorithm for obtaining the maximum constrained supereigenvalue and depicting the feasible region of all the constrained supereigenvectors for a given matrix. It also gives some examples for illustrating the algorithms, respectively.

\end{quote}

\section{Introduction}
An eigenproblem is a very classical and important research topic in linear algebra, both in theory and in practice. The eigenproblem in linear algebra is to find the eigenvalues and the eigenvectors of a given matrix $A$. A scalar $\lambda$ is called an eigenvalue of $A$, if there exists a nonzero vector $x$ satisfying
 \begin{equation}\label{eq1} Ax^T = \lambda x^T,\end{equation}
where the operations of equation (\ref{eq1}) are the common plus (+) and multiplication ($\times$) in the real number set.
In some practical applications, the eigenproblem is used to describe the steady states of discrete event systems, for instance, the eigenproblems in fuzzy algebra and max-plus algebra \cite{co,ol,lu}. Now, the eigenproblem in algebra has been studied by many authors, such as max-plus algebra \cite{bu,ba,ce05,cu62,cu79,ga16}, max-min algebra \cite{ce92,ce97,cu79, cu91,ga02,ga014,sa}, max-prod algebra \cite{ra016,el,bi,gu,pe}, max-Lukasiewicz and max-drastic algebra \cite{ra12,ga14}.

In particular, Cuninghame-Green \cite{cu62} first studied the eigenproblem in max-plus algebra, and he even gave all solutions of the eigenproblem with irreducible matrix \cite{cu79}. Bapat et al. \cite{ba} further proved the spectral theorem for the reducible matrices. In \cite{ce05}, Cechlrov defined the universal and possible eigenvectors and proposed an algorithm for judging whether a given vector is a possible eigenvector and a universal eigenvector exists. Gavalec et al. \cite{ga16} investigated the eigenproblem of a given matrix with interval coefficients and supplied polynomial algorithms to identify three types of tolerance interval eigenvectors.

In max-min algebra, an algorithm for calculating the maximum eigenvector of a given matrix is proposed \cite{ce92,ce97,cu79, cu91}. Gavalec showed that the complete structure of the eigenspace is a union of permutation non-decreasing eigenvector intervals \cite{ga02}. In max-drastic algebra, Gavalec et al. \cite{ga14} completely described the structure of the eigenspace of a given matrix. Rashid et al. \cite{ra12} obtained a similar result for a square matrix in max-{\L}ukasiewicz algebra. Rashid et al. \cite{ra016} also investigated the eigenspace of a given matrix of order $3$ in max-prod algebra.

The theory of fuzzy relation inequalities has been widely used in BitTorrent-like peer-to-peer (P2P) file sharing system. Assume that the P2P file sharing system has $n$ terminals, which are denoted by $A_1, A_2, \cdots, A_n$. Each terminal should share its local file resources to any other terminal.  At the same time, the file data can be downloaded from any other terminal. Suppose that the $j$th terminal $A_j$ sends the file data to another terminal with quality level $x_j$ and the bandwidth between $A_i$ and $A_j$ is $a_{ij}$. For data transmission, the download quality level of $A_i$ from $A_j$ is actually $a_{ij}\wedge x_j$ due to the bandwidth limitation, on which the file data of terminal $A_i$ from other terminals is $a_{i1}\wedge x_1 + \cdots +a_{i{}i-1}\wedge x_{i-1}+a_{i{i+1}}\wedge x_{i+1}+\cdots+ a_{in}\wedge x_n$. In order to satisfy the downloading quality requirements of users, total download quality of $A_i$ should be no less than $b_i$ ($b_i > 0$). Then the P2P file sharing system is reduced to a system of fuzzy relation inequations with addition-min composition as follows.

\begin{equation}\label{b}
\left\{\begin{array}{ll}
a_{11}\wedge x_1 + a_{12}\wedge x_{2} + \cdots + a_{1n}\wedge x_n\geq b_1,\\
a_{21}\wedge x_1 + a_{22}\wedge x_{2} + \cdots + a_{2n}\wedge x_n\geq b_2,\\
\cdots\\
a_{n1}\wedge x_1 + a_{n2}\wedge x_{2} + \cdots + a_{nn}\wedge x_n\geq b_n,
\end{array}\right.
\end{equation}
 where $a_{ij},x_{j} \in [0,1]$, $b_i > 0$ $(i = 1, 2, \cdots, n; j = 1, 2, \cdots, n)$, $a_{ij} \wedge x_{j} = \min \{a_{ij}, x_j\}$, and the operation $'+'$ is the ordinary addition \cite{Li2012}. System (\ref{b}) can be tersely described as follows
\begin{equation*}A \odot {x^T} \ge {b^T}\end{equation*}
where $A=(a_{ij})_{n \times n}$, $x=(x_1, x_2, \cdots, x_n)$, $b=(b_1, b_2, \cdots, b_n)$ and
$(a_{i1}, a_{i2}, \cdots , a_{in}) \odot ( x_1, x_2, \cdots, x_n) = a_{i1} \wedge x_1 + a_{i2} \wedge x_2 + \cdots + a_{in} \wedge x_n$.

In order to avoid network congestion and improve the stability of data transmission, we need to enhance the relevance between the data download quality and the data sending quality of the terminal. Following this idea, a popular method is to give the data download quality a proportional to the data sending quality \cite{Yang2019}. This means that the data download quality to the data sending quality will reach a steady state. We say that the system reaches a steady regime. If the proportionality coefficient is denoted by $\lambda\geq 0$, then the corresponding steady regime of the P2P file sharing system can be written in mathematics as follows.
\begin{equation}\label{x}
\left\{\begin{array}{ll}
a_{11}\wedge x_1 + a_{12}\wedge x_{2} + \cdots + a_{1n}\wedge x_n=\lambda x_1,\\
a_{21}\wedge x_1 + a_{22}\wedge x_{2} + \cdots + a_{2n}\wedge x_n=\lambda x_2,\\
\cdots\\
a_{n1}\wedge x_1 + a_{n2}\wedge x_{2} + \cdots + a_{nn}\wedge x_n=\lambda x_n.
\end{array}\right.
\end{equation}
The matrix form of system (\ref{x}) is
$$A \odot x^T = \lambda x^T$$
 where $\lambda\in [0,+\infty)$. Considering the addition-min composition $\odot$, $\lambda x^T$ means that $\lambda$ multiplies by $x_j$, $j\in\{1,2 \cdots, n\}$. It is clear that the steady states of the data download quality to the data sending quality are more favorable to the terminal. In other words, a steady value will benefit users. This article aims to obtain the steady value $\lambda$ and the steady solution $x$, called an eigenvalue and eigenvector of $A$, respectively.

 The rest of this article is organized as follows. In Section \ref{se2}, we present some necessary notation and known results for following this article. In Section \ref{se3}, we introduce an eigenproblem of addition-min algebra, provide a sufficient and necessary condition for a vector being an eigenvector of a given matrix $A$, and then present an algorithm for finding all the eigenvalues and eigenvectors of a given matrix $A$ for system (\ref{x}). In Section \ref{se4}, we consider the so-called constrained eigenproblem, show a sufficient and necessary condition for a vector being a constrained eigenvector of a given matrix $A$, and give an algorithm for computing all the constrained eigenvectors and eigenvalues of a given matrix $A$. In Section 5, we make the algorithms obtained in Sections 3 and 4 be suitable for the supereigenvectors and constrained supereigenvectors of addition-min algebras, respectively. A concluding remark is drawn in Section \ref{se6}.

 \section{Preliminaries}\label{se2}
This section presents some basic notation and known results.

We call the algebra $([0,1], +, \wedge)$ an addition-min algebra and denote the real unit interval by $\mathcal{I}= [0, 1]$. For a given natural number $n$, we write $N= \{1, 2,\cdots,n\}$. Further, the notation $\mathcal{I}^{n\times n}$ denotes the set of $n\times n$ matrices over $\mathcal{I}$, the set of $1\times n$ vectors over $\mathcal{I}$ is denoted by $\mathcal{I}^n$ and $\theta=(0,0, \cdots, 0)$.

For $x=(x_1, x_2, \cdots, x_n), y=(y_1, y_2, \cdots, y_n)\in \mathcal{I}^n$, define $x\leq y$ if and only if $x_i\leq y_i$ for arbitrary $i\in N$, and define $x < y$ if and only if $x_i \leq y_i$ for arbitrary $i\in N$ and there is a $j\in N$ such that $x_j < y_j$.
\begin{lemma}[\cite{Li2012,Yang2014}]\label{Le2.1}
System (\ref{b}) is solvable if and only if $(1, 1,\cdots, 1)$ is a solution of system (\ref{b}).
 \end{lemma}
\begin{lemma}[\cite{Li2012,Yang2014}]\label{lem2}
Let $x = (x_1, x_2, \cdots, x_n)$ be  a solution of system (\ref{b}). Then we have:
\begin{enumerate}
\item [(1)] $x > 0$.
\item [(2)] For any $i \in N$, $j \in N$, $x_j \ge b_i - \sum\limits_{k \in N \setminus \{j\}} {a_{ik} \wedge x_k} \ge  b_i - \sum\limits_{k \in N \setminus \{ j \}} {a_{ik}} $.
\item [(3)] For any $i \in N$, $j \in N$, $a_{ij} \ge b_i - \sum\limits_{k \in N \setminus \{ j \}} {a_{ik} \wedge x_k} \ge b_i - \sum\limits_{k \in N \setminus \{ j\}} a_{ik}$.
\end{enumerate}
\end{lemma}

Let $\check{\alpha}= (\check{\alpha}_1, \check{\alpha}_2, \cdots, \check{\alpha}_n)$ with
\begin{equation*}
 \check{\alpha}_j = \max \{0, b_i-\sum \limits_{k \in N \setminus \{ j \}}a_{ik}|i\in N\}
\end{equation*}and $\hat{\alpha} = (\hat{\alpha}_1, \hat{\alpha}_2, \cdots, \hat{\alpha}_n)$ with
\begin{equation*}\label{eq3}
 \hat{\alpha}_j= \max \{a_{ij}|i\in N\}.
\end{equation*}
Then we have the following lemma.
\begin{lemma}[\cite{Yang2014}]\label{Le2.2}
 For any solution $x$ of system (\ref{b}), it holds that $\check{\alpha}\leq x$.
\end{lemma}

\section{Eigenproblems of addition-min algebras}\label{se3}
In this section, we first introduce an eigenproblem of addition-min algebra, and then investigate the conditions for a vector being an eigenvector of a given matrix $A$. We also explore an algorithm for finding all the eigenvalues and eigenvectors of a given matrix $A$ for system (\ref{x}).

In addition-min algebra, for a given matrix $A\in \mathcal{I}^{n\times n}$, the task of finding a vector $x\in \mathcal{I}^n$ with $x\neq \theta$ and a scalar $\lambda\in [0,+\infty)$ satisfying system \eqref{x} is called an eigenproblem of addition-min algebra. The scalar $\lambda$ is called an eigenvalue of $A$, and the corresponding vector $x$ is called an eigenvector of $A$ associated with $\lambda$.

The notation $V(A,\lambda)$ denotes the set consisting of all eigenvectors of $A$ corresponding to $\lambda$, and $\Lambda(A)$ denotes the set consisting of all eigenvalues of $A$, i.e., $$V(A,\lambda)= \{x\in \mathcal{I}^n|A \odot x^T= \lambda x^T\mbox{ and } x\neq \theta\}$$
and $$\Lambda(A)=\{\lambda\in [0,+\infty)|V(A,\lambda)\}.$$
Let $|\emptyset|=0$ and denote the set of all eigenvectors of $A$ by $V(A)$, i.e., $$V(A)=\bigcup\limits_{\lambda\in \Lambda(A)}V(A,\lambda).$$

In what follows, we discuss the eigenproblems of addition-min algebras.
\begin{definition}\label{de3.1}
For $A=(a_{ij})_{n\times n}\in \mathcal{I}^{n\times n}$ and $j\in N$, denote $|\{a_{ij}| 0< a_{ij}<1, i\in N\}|=t_j$ and $Q_j = \{q_{0j}, q_{1j}, q_{2j}, \cdots, q_{t_jj}, q_{(t_j+1)j}\}$ with $0= q_{0j}<  q_{1j} <  q_{2j}< \cdots < q_{t_jj}<q_{(t_j+1)j}=1$, where $q_{kj}\in \{a_{ij}| i\in N\}, k= 1, 2, \cdots, t_j$.
\end{definition}

 Let $K = \{(k_1, k_2, \cdots, k_n) |k_j\in K_j \mbox{ for any }j\in N\}$ with
\begin{equation}\label{k}
   K_j= \{k|0< q_{kj}\mbox{ and } q_{kj}\in Q_j\}.
\end{equation}

 Notice that from Definition \ref{de3.1}, one can check that $K_j\neq \emptyset$ and $|Q_j|\geq 2$ for any $j\in N$. In particular, if $t_j=0$ for a $j\in N$ then $Q_j = \{0,1\}$.
 \begin{theorem}\label{th3.1}
 Let $x = (x_1, x_2, \cdots, x_n)\in \mathcal{I}^n$. Then $x \in V(A, \lambda)$ if and only if there exists a $k=(k_1, k_2, \cdots, k_n)\in K$ such that
\begin{equation}\label{ab}\left\{\begin{array}{ll}
q_{{(k_j-1)}j} \leq x_j \leq q_{k_jj}\mbox{ for all } j\in N,\\ \\
\sum \limits_{j\in N}[\delta_{ij}\cdot x_j + (1 - \delta_{ij})\cdot a_{ij}]= \lambda x_i\mbox{ for all } i\in N
\end{array}\right.
\end{equation}
where the operation $``\cdot"$ represents the ordinary multiplication and $$\delta_{ij} = \left\{\begin{array}{l}
1, q_{k_jj} \leq a_{ij},\\
0, a_{ij} \leq q_{(k_j-1)j}.
\end{array} \right.$$
\end{theorem}
\begin{proof}
Suppose that $x = (x_1, x_2, \cdots, x_n)\in V(A, \lambda)$. According to Definition \ref{de3.1} and Formula (\ref{k}), it is clear that there is a $k_j$ such that $q_{{(k_j-1)}j} \leq x_j \leq q_{k_jj}$ for any $j\in N$. Thus, there exists a $k=(k_1, k_2, \cdots, k_n)\in K$ such that $q_{{(k_j-1)}j} \leq x_j \leq q_{k_jj}$ for any $j\in N$. For any $i\in N$, we get
\begin{eqnarray*}
\lambda x_i&=&\sum \limits_{j\in N} a_{ij}\wedge x_j
\\&=&\sum \limits_{j\in N, q_{k_jj} \leq a_{ij}}a_{ij}\wedge x_j+ \sum \limits_{j\in N, a_{ij}\leq q_{(k_j-1)j}}a_{ij}\wedge x_j
\\&=&\sum \limits_{j\in N, q_{k_jj} \leq a_{ij}} x_j+ \sum \limits_{j\in N, a_{ij}\leq q_{(k_j-1)j}}a_{ij}
\\ &=&\sum \limits_{j\in N, q_{k_jj} \leq a_{ij}}\delta_{ij}\cdot x_j + \sum \limits_{j\in N, a_{ij}\leq q_{(k_j-1)j}}(1 - \delta_{ij})\cdot a_{ij}
\\&=&\sum \limits_{j\in N}[\delta_{ij}\cdot x_j + (1 - \delta_{ij})\cdot a_{ij}],
 \end{eqnarray*}
where the operation $``\cdot"$ represents the ordinary multiplication and $$\delta_{ij} = \left\{ \begin{array}{l}
1, q_{k_jj} \leq a_{ij},\\
0, a_{ij} \leq q_{(k_j-1)j}.
\end{array} \right.$$

Now, suppose that there exists a $k=(k_1, k_2, \cdots, k_n)\in K$ such that $x$ satisfies system (\ref{ab}). Then for any $i\in N$, we have
\begin{eqnarray*}
\lambda x_i&=&\sum \limits_{j\in N}[\delta_{ij}\cdot x_j + (1 - \delta_{ij})\cdot a_{ij}]
\\ &=&\sum \limits_{j\in N, q_{k_jj} \leq a_{ij}}\delta_{ij}\cdot x_j + \sum \limits_{j\in N, a_{ij}\leq q_{(k_j-1)j}}(1 -\delta_{ij})\cdot a_{ij}
\\&=&\sum \limits_{j\in N, q_{k_jj} \leq a_{ij}} x_j+ \sum \limits_{j\in N, a_{ij}\leq q_{(k_j-1)j}}a_{ij}
\\&=&\sum \limits_{j\in N, q_{k_jj} \leq a_{ij}}a_{ij}\wedge x_j+ \sum \limits_{j\in N, a_{ij}\leq q_{(k_j-1)j}}a_{ij}\wedge x_j
\\&=&\sum \limits_{j\in N} a_{ij}\wedge x_j.
 \end{eqnarray*}
 Therefore, $x\in V(A, \lambda)$.
\end{proof}

Note that in Theorem \ref{th3.1}, $\lambda$ can be seen as a given constant with $\lambda\in [0, +\infty)$. So that every solution $x= (x_1, x_2, \cdots, x_n)$ of system (\ref{ab}) can be represented by $\lambda$, and the corresponding $\lambda$ can be determined by the $x= (x_1, x_2, \cdots, x_n)$ since $q_{{(k_j-1)}j} \leq x_j \leq q_{k_jj}\mbox{ for all } j\in N$.

Furthermore, Theorem \ref{th3.1} implies the following two statements.
\begin{theorem}\label{th3.2}
If system (\ref{ab}) is unsolvable for any $k=(k_1, k_2, \cdots, k_n)\in K$ then $V(A, \lambda)=\emptyset$.
\end{theorem}
\begin{theorem}\label{th3.2}
For any $k=(k_1, k_2, \cdots, k_n)\in K$, if $x= (x_1, x_2, \cdots, x_n)$ satisfies system (\ref{ab}), then $x\in V(A, \lambda)$.
\end{theorem}

For any $k=(k_1, k_2, \cdots, k_n)\in K$, we use $V(A, \lambda, k)$ denotes the solution set of system (\ref{ab}) corresponding to $k=(k_1, k_2, \cdots, k_n)\in K$. Then based on Theorems \ref{th3.1} and \ref{th3.2}, the eigenproblems of addition-min algebras are equivalent to solving system (\ref{ab}). Therefore, we can summarize an algorithm for finding all the eigenvectors and eigenvalues of a given $A$ as follows.
\begin{algorithm}\label{al1}
Input $A=(a_{ij})_{n\times n}$. Output $\Lambda(A)$ and $V(A)$.

Step 1. Compute $K = \{(k_1, k_2, \cdots, k_n) |k_j\in K_j \mbox{ for any }j\in N\}$ defined by (\ref{k}).

Step 2. For any $k=(k_1, k_2, \cdots, k_n)\in K$, construct the corresponding system (\ref{ab}).

Step 3. Compute $V(A, \lambda, k)$ by the corresponding system (\ref{ab}).

Step 4. Output $\Lambda(A)=\bigcup \limits_{k\in K}\{\lambda\in [0,+\infty)|V(A,\lambda, k)\}$ and $V(A)=\bigcup\limits_{\lambda\in \Lambda(A)}V(A,\lambda)$.

Step 5. End.
\end{algorithm}
\begin{theorem}\label{th3.3} Algorithm \ref{al1} terminates after $O((n+1)^{n})$ operations.
\end{theorem}
\begin{proof} In Step 1, it costs $\sum\limits_{j\in N}(t_j+2)$ operations for computing $K$, where $\sum\limits_{j\in N}(t_j+2)\leq (n+2)n$. In Steps 2 and 3, for any $k\in K $, it takes $6n^2+3n$ operations for solving a corresponding system (\ref{ab}). Therefore, it costs $(6n^2+3n)|K|$ operations, where $|K|=\prod\limits_{j\in N}(t_j+1)\leq (n+1)^n$. Step 4 needs $2(|K|-1)$ operations for computing $\Lambda(A)$ and $V(A)$. Therefore, the amount of computation of Algorithm \ref{al1} is
$\sum\limits_{j\in N}(t_j+2)+(6n^2+3n)|K|+2(|K|-1)=(6n^2+3n+2)|K|+\sum\limits_{j\in N}(t_j+2)-2 $. Consequently, the total computational complexity of Algorithm \ref{al1} is $O((n+1)^{n})$.\end{proof}

The following example illustrates Algorithm \ref{al1}.
\begin{example}\label{ex1}
\emph{Consider the following system:}
$$\left\{ \begin{array}{l}
0.4\wedge x_1 + 0.6\wedge x_2=\lambda x_1,\\
0.2\wedge x_1 + 0.5\wedge x_2=\lambda x_2.\\
\end{array} \right.$$
\end{example}
Step 1. Compute $K_1=\{1, 2, 3\}$ and $K_2=\{1, 2, 3\}$. Then
\begin{eqnarray*}
K=\{(1, 1), (1, 2),(1, 3), (2, 1), (2, 2), (2, 3), (3, 1), (3, 2), (3, 3)\}.
\end{eqnarray*}
Step 2. Construct and solve the corresponding system (\ref{ab}):
\begin{enumerate}
\item [(1)] For $k_1=(1, 1)$,
$$\left\{ \begin{array}{l}
x_1 + x_2 =\lambda x_1,\\
x_1 + x_2=\lambda x_2,\\
0\leq x_1\leq 0.2,\\
0\leq x_2\leq 0.5.
\end{array} \right.$$
We get $x_1=x_2$ and $\lambda=2$. Let $x_1=x_2=t$. Since $0 \leq x_1\leq 0.2$ and $0 \leq x_2\leq 0.5$, we have $t\in [0, 0.2]$. When $t=0$, $x=(0, 0)\notin V(A, \lambda)$. Therefore, $V(A, \lambda, k_1)=\{(t,t)|t\in (0, 0.2]\mbox{ and } \lambda=2\}$.

\item [(2)] For $k_2=(1, 2)$,
$$\left\{ \begin{array}{l}
x_1 + x_2 =\lambda x_1,\\
x_1 +0.5=\lambda x_2,\\
0\leq x_1\leq 0.2,\\
0.5\leq x_2\leq 0.6.
\end{array} \right.$$
We get $x_1=\frac{0.5}{\lambda^2-\lambda-1}$, $x_2=\frac{0.5(\lambda-1)}{\lambda^2-\lambda-1}$ and $\lambda^2-\lambda-1>0$. Since $0 \leq \frac{0.5}{\lambda^2-\lambda-1}\leq 0.2$ and $0.5\leq \frac{0.5(\lambda-1)}{\lambda^2-\lambda-1}\leq 0.6$, it is impossible. Therefore, $V(A, \lambda, k_2)=\emptyset$.

\item [(3)] For $k_3=(1, 3)$,
$$\left\{ \begin{array}{l}
x_1 + 0.6 =\lambda x_1,\\
x_1 +0.5=\lambda x_2,\\
0\leq x_1\leq 0.2,\\
0.6\leq x_2\leq1.
\end{array} \right.$$
We get $x_1=\frac{0.6}{\lambda-1}$, $x_2=\frac{0.5\lambda+0.1}{\lambda(\lambda-1)}$ and $\lambda-1>0$. Since $0\leq \frac{0.6}{\lambda-1}\leq 0.2$ and $0.6\leq \frac{0.5\lambda+0.1}{\lambda(\lambda-1)}\leq 1$, it is impossible. Therefore, $V(A, \lambda, k_3)=\emptyset$.

\item [(4)] For $k_4=(2, 1)$,
$$\left\{ \begin{array}{l}
x_1 + x_2 =\lambda x_1,\\
0.2 + x_2=\lambda x_2,\\
0.2\leq x_1\leq 0.4,\\
0\leq x_2\leq 0.5.
\end{array} \right.$$
We get $x_1=\frac{0.2}{(\lambda-1)^2}$, $x_2=\frac{0.2}{\lambda-1}$ and $\lambda-1>0$. Since $0.2\leq \frac{0.2}{(\lambda-1)^2}\leq 0.4$ and $0\leq \frac{0.2}{\lambda-1}\leq 0.5$, we have $\frac{\sqrt{2}+2}{2}\leq\lambda\leq 2$. Therefore, $V(A, \lambda, k_4)=\{(\frac{0.2}{(\lambda-1)^2}, \frac{0.2}{\lambda-1})|\frac{\sqrt{2}+2}{2}\leq\lambda\leq 2\}.$

\item [(5)] For $k_5=(2, 2)$,
$$\left\{ \begin{array}{l}
x_1 + x_2 =\lambda x_1,\\
0.2 + 0.5=\lambda x_2,\\
0.2\leq x_1\leq 0.4,\\
0.5\leq x_2\leq 0.6.
\end{array} \right.$$
We get $x_1=\frac{0.7}{\lambda(\lambda-1)}$, $x_2=\frac{0.7}{\lambda}$ and $\lambda-1>0$. Since $0.2\leq \frac{0.7}{\lambda(\lambda-1)}\leq 0.4$ and $0.5\leq \frac{0.7}{\lambda}\leq 0.6$, it is impossible. Therefore, $V(A, \lambda, k_5)=\emptyset.$

\item [(6)]  For $k_6=(2, 3)$,
$$\left\{ \begin{array}{l}
x_1 + 0.6 =\lambda x_1,\\
0.2 + 0.5=\lambda x_2,\\
0.2\leq x_1\leq 0.4,\\
0.6\leq x_2\leq 1.
\end{array} \right.$$
We get $x_1=\frac{0.6}{\lambda-1}$, $x_2=\frac{0.7}{\lambda}$ and $ \lambda-1>0$. Since $0.2\leq \frac{0.6}{\lambda-1}\leq 0.4$ and $0.6\leq \frac{0.7}{\lambda}\leq 1$, it is impossible. Therefore, $V(A, \lambda, k_6)=\emptyset.$

\item [(7)] For $k_7=(3, 1)$,
$$\left\{ \begin{array}{l}
0.4 + x_2 =\lambda x_1,\\
0.2 + x_2=\lambda x_2,\\
0.4\leq x_1\leq 1,\\
0\leq x_2\leq 0.5.
\end{array} \right.$$
We get $x_1=\frac{0.4\lambda-0.2}{\lambda(\lambda-1)}$, $x_2=\frac{0.2}{\lambda-1}$ and $\lambda-1>0$. Since $0.4\leq \frac{0.4\lambda-0.2}{\lambda(\lambda-1)}\leq 1$ and $0\leq \frac{0.2}{\lambda-1}\leq 0.5$, we have $\frac{7}{5}\leq\lambda\leq \frac{\sqrt{2}+2}{2}$. Therefore, $V(A, \lambda, k_7)=\{(\frac{0.4\lambda-0.2}{\lambda(\lambda-1)}, \frac{0.2}{\lambda-1})|\frac{7}{5}\leq\lambda\leq \frac{\sqrt{2}+2}{2}\}.$

\item [(8)] For $k_8=(3, 2)$,
$$\left\{ \begin{array}{l}
0.4 + x_2 =\lambda x_1,\\
0.2 + 0.5=\lambda x_2,\\
0.4\leq x_1\leq 1,\\
0.5\leq x_2\leq 0.6.
\end{array} \right.$$
We get $x_1=\frac{0.4\lambda+0.7}{\lambda^2}$, $x_2=\frac{0.7}{\lambda}$ and $\lambda>0$. Since $0.4\leq \frac{0.4\lambda+0.7}{\lambda^2}\leq 1$ and $0.5\leq \frac{0.7}{\lambda}\leq 0.6$, we have $\frac{7}{6}\leq\lambda\leq \frac{7}{5}$. Therefore, $V(A, \lambda, k_8)=\{(\frac{0.4\lambda+0.7}{\lambda^2}, \frac{0.7}{\lambda})|\frac{7}{6}\leq\lambda\leq \frac{7}{5})\}.$

\item [(9)] For $k_9=(3, 3)$,
$$\left\{ \begin{array}{l}
0.4 + 0.6 =\lambda x_1,\\
0.2 + 0.5=\lambda x_2,\\
0.4\leq x_1\leq 1,\\
0.6\leq x_2\leq 1.
\end{array} \right.$$
We get $x_1=\frac{1}{\lambda}$, $x_2=\frac{0.7}{\lambda}$ and $\lambda>0$. Since $0.4\leq \frac{1}{\lambda}\leq 1$ and $0.6\leq \frac{0.7}{\lambda}\leq 1$, we have $1\leq\lambda\leq \frac{7}{6}$. Therefore, $V(A, \lambda, k_9)=\{(\frac{1}{\lambda}, \frac{0.7}{\lambda})|1\leq\lambda\leq \frac{7}{6}\}.$

\end{enumerate}
Step 3. Output
\begin{eqnarray*}
\Lambda(A)&=&\bigcup \limits_{k\in K}\{\lambda\in [0,+\infty)|V(A,\lambda, k)\}
\\&=&\{\lambda\in [0,+\infty)|V(A,\lambda, k_1)\}\cup \{\lambda\in [0,+\infty)|V(A,\lambda, k_4)\}\cup \{\lambda\in [0,+\infty)|V(A,\lambda, k_7)\}\\&\quad&\cup \{\lambda\in [0,+\infty)|V(A,\lambda, k_8)\}\cup \{\lambda\in [0,+\infty)|V(A,\lambda, k_9)\}
\\&=&\{\lambda|\lambda=2\}\cup\{\lambda|\frac{\sqrt{2}+2}{2}\leq\lambda\leq 2\}\cup\{\lambda|\frac{7}{5}\leq\lambda\leq \frac{\sqrt{2}+2}{2}\}\cup\{\lambda|\frac{7}{6}\leq\lambda\leq \frac{7}{5}\}\\&&\cup\{\lambda|1\leq\lambda\leq \frac{7}{6}\}
\\&=& \{\lambda|1\leq\lambda\leq 2\}
\end{eqnarray*}
and
\begin{eqnarray*}
 V(A)&=&\bigcup\limits_{\lambda\in \Lambda(A)}V(A,\lambda)
 \\&=&\{(t,t)|t\in (0, 0.2]\mbox{ and }\lambda=2\}\cup \{(\frac{0.2}{(\lambda-1)^2}, \frac{0.2}{\lambda-1})|\frac{\sqrt{2}+2}{2}\leq\lambda\leq 2\}\\&&\cup\{(\frac{0.4\lambda-0.2}{\lambda(\lambda-1)}, \frac{0.2}{\lambda-1})|\frac{7}{5}\leq\lambda\leq \frac{\sqrt{2}+2}{2}\}\cup \{(\frac{0.4\lambda+0.7}{\lambda^2}, \frac{0.7}{\lambda})|\frac{7}{6}\leq\lambda\leq \frac{7}{5}\}\\&&\cup \{(\frac{1}{\lambda}, \frac{0.7}{\lambda})|1\leq\lambda\leq \frac{7}{6}\}.
\end{eqnarray*}
\section{Constrained eigenproblems of addition-min algebras}\label{se4}
In order to satisfy the downloading quality requirements of users and improve the stability of data transmission in the P2P file sharing system, it is worth noting the steady states of the data download quality to the data sending quality under system (\ref{b}). Based on such a consideration, we investigate the eigenproblem of matrix in addition-min algebra under system (\ref{b}) in this section.

For a given $A\in \mathcal{I}^{n\times n}$, the task of finding a vector $x\in \mathcal{I}^n$ with $x\neq \theta$ and a scalar $\lambda$ satisfying both systems \eqref{b} and \eqref{x} is called a constrained eigenproblem of addition-min algebra. The scalar $\lambda$ is called a constrained eigenvalue of $A$, and the corresponding vector $x$ is called a constrained eigenvector of $A$ associated with $\lambda$.

Denoted the set consisting of all constrained eigenvectors of $A$ corresponding to $\lambda$ by $V^*(A,\lambda)$, and the set consisting of all constrained eigenvalues of $A$ by $\Lambda^*(A)$, i.e.,
 $$V^*(A,\lambda)= \{x\in \mathcal{I}^n|A \odot x^T \geq b^T, A \odot x^T= \lambda x^T\mbox{ and } x\neq \theta\}$$
and $$\Lambda^*(A)=\{\lambda\in [0,+\infty)|V^*(A,\lambda)\}.$$
We also denote by $V^*(A)$ the set consisting of all constrained eigenvectors of $A$, i.e., $$V^*(A)=\bigcup\limits_{\lambda\in \Lambda^*(A)}V^*(A,\lambda).$$

From Lemma \ref{lem2}, we have $\check{\alpha}\leq \hat{\alpha}$. Therefore, for any $j\in N$ there exists an $i\in N$ such that $\check{\alpha}_j\leq a_{ij}$. In this way, we can give the following definition.
\begin{definition}\label{de4.1}
For $A=(a_{ij})_{n\times n}\in \mathcal{I}^{n\times n}$ and $j\in N$, denote $D_j = \{d_{0j}, d_{1j}, d_{2j}, \cdots, d_{l_jj}\}$, where $D_j$ satisfies the following two conditions:
 \begin{enumerate}
\item [(i)] $\check{\alpha}_j= d_{0j}< d_{1j} < d_{2j}< \cdots < d_{l_jj}=1$, where $d_{pj}\in \{a_{ij}| i\in N\}, p= 1, 2, \cdots, l_j-1$.
\item [(ii)] For any $i\in N$, if $\check{\alpha}_j\leq a_{ij}$ then there exists a unique $p\in \{0, 1, 2, \cdots, l_j\}$ such that $a_{ij}=d_{pj}$.
\end{enumerate}
\end{definition}

Denote
 \begin{equation*}\label{eqj}
 N^*=\{j\in N|\check {\alpha}_j=1\}
  \end{equation*}
and $P = \{(p_1, p_2, \cdots, p_n) |p_j\in P_j \mbox{ for any }j\in N\}$ with
\begin{equation}\label{p}
  P_j= \left\{\begin{array}{ll}
\{0\}, j\in N^*,\\ \\
\{p|\check{\alpha}_j< d_{pj}\mbox{ and } d_{pj}\in D_j\}, j\in N\setminus N^*.
\end{array}\right.
\end{equation}

Then we have the following theorem.
 \begin{theorem}\label{th4.1}
 Let $x = (x_1, x_2, \cdots, x_n)\in \mathcal{I}^n$. Then $x \in V^*(A, \lambda)$ if and only if there exists a $p=(p_1, p_2, \cdots, p_n)\in P$ such that
\begin{equation}\label{ac}\left\{\begin{array}{ll}
x_j= d_{0j}=1\mbox{ for all } j\in N^*,\\\\
d_{{(p_j-1)}j} \leq x_j \leq d_{p_jj}\mbox{ for all } j\in N\setminus N^*,\\ \\
\sum \limits_{j\in N^*}a_{ij}+\sum \limits_{j\in N\setminus N^*}[\gamma_{ij}\cdot x_j + (1 - \gamma_{ij})\cdot a_{ij}]\geq b_i\mbox{ for all } i\in N,\\\\
\sum \limits_{j\in N^*}a_{ij}+\sum \limits_{j\in N\setminus N^*}[\gamma_{ij}\cdot x_j + (1 - \gamma_{ij})\cdot a_{ij}]= \lambda x_i\mbox{ for all } i\in N
\end{array}\right.
\end{equation}
where the operation $``\cdot"$ represents the ordinary multiplication and
$$\gamma_{ij} = \left\{ \begin{array}{l}
1, d_{p_jj} \leq a_{ij},\\
0, a_{ij} \leq d_{(p_j-1)j}.
\end{array} \right.$$
\end{theorem}
\begin{proof}
Suppose that $x = (x_1, x_2, \cdots, x_n)\in V^*(A, \lambda)$. From Lemma \ref{Le2.2}, $\check {\alpha}_j\leq  x_j\leq 1$ for any $j\in N$, then $\check {\alpha}_j=1$ or $\check {\alpha}_j<1$. If $\check {\alpha}_j=1$, then $x_j=1$ and $j\in N^*$. According to Definition \ref{de4.1} and Formula (\ref{p}), $p_j=0\in P_j$ and $\check {\alpha}_j=x_j=d_{0j}=1$ for any $j\in N^*$. If $\check {\alpha}_j<1$ then $j\in N\setminus N^*$. Thus by Definition \ref{de4.1} and Formula (\ref{p}), $\check {\alpha}_j\leq d_{(p_j-1)j}\leq x_j\leq d_{p_jj}\leq1$, i.e., there exists a $p_j\in P_j$ such that $d_{(p_j-1)j}\leq x_j\leq d_{p_jj}$ for any $j\in N\setminus N^*$. Therefore, $P_j\neq \emptyset$ for any $j\in N$ and
\begin{eqnarray*}
\sum \limits_{j\in N} a_{ij}\wedge x_j&=&\sum \limits_{j\in N^*} a_{ij}\wedge x_j+\sum \limits_{j\in N\setminus N^*} a_{ij}\wedge x_j
\\&=&\sum \limits_{j\in N^*} a_{ij}\wedge x_j+\sum \limits_{j\in N\setminus N^*, d_{p_jj} \leq a_{ij}}a_{ij}\wedge x_j+ \sum \limits_{j\in N\setminus N^*, a_{ij}\leq d_{(p_j-1)j}}a_{ij}\wedge x_j
\\&=&\sum \limits_{j\in N^*} a_{ij}+\sum\limits_{j\in N\setminus N^*, d_{p_jj} \leq a_{ij}} x_j+ \sum \limits_{j\in N\setminus N^*, a_{ij}\leq d_{(p_j-1)j}}a_{ij}
\\ &=&\sum \limits_{j\in N^*} a_{ij}+\sum \limits_{j\in N\setminus N^*, d_{p_jj} \leq a_{ij}}\gamma_{ij}\cdot x_j + \sum \limits_{j\in N\setminus N^*, a_{ij}\leq d_{(p_j-1)j}}(1 - \gamma_{ij})\cdot a_{ij}
\\&=&\sum \limits_{j\in N^*} a_{ij}+\sum \limits_{j\in N\setminus N^*}[\gamma_{ij}\cdot x_j + (1 - \gamma_{ij})\cdot a_{ij}],
 \end{eqnarray*}
where the operation $``\cdot"$ represents the ordinary multiplication and
 $$\gamma_{ij} = \left\{ \begin{array}{l}
1, d_{p_jj} \leq a_{ij},\\
0, a_{ij} \leq d_{(p_j-1)j}.
\end{array} \right.$$

Because of $x = (x_1, x_2, \cdots, x_n)\in V^*(A, \lambda)$, we have $b_i\leq \sum \limits_{j\in N} a_{ij}\wedge x_j=\sum \limits_{j\in N^*} a_{ij}+\sum \limits_{j\in N\setminus N^*}[\gamma_{ij}\cdot x_j + (1 - \gamma_{ij})\cdot a_{ij}]$ and $\lambda x_i= \sum \limits_{j\in N} a_{ij}\wedge x_j=\sum \limits_{j\in N^*} a_{ij}+\sum \limits_{j\in N\setminus N^*}[\gamma_{ij}\cdot x_j + (1 - \gamma_{ij})\cdot a_{ij}]$ for any $i\in N$, i.e., $\sum \limits_{j\in N^*}a_{ij}+\sum \limits_{j\in N\setminus N^*}[\gamma_{ij}\cdot x_j + (1 - \gamma_{ij})\cdot a_{ij}]\geq b_i$ and $\sum \limits_{j\in N^*}a_{ij}+\sum \limits_{j\in N\setminus N^*}[\gamma_{ij}\cdot x_j + (1 - \gamma_{ij})\cdot a_{ij}]= \lambda x_i$ for all $i\in N$.

Conversely, suppose that there exists a $p=(p_1, p_2, \cdots, p_n)\in P$ such that $x$ satisfies system (\ref{ac}). Then
\begin{eqnarray*}
&\quad&\sum \limits_{j\in N^*} a_{ij}+\sum \limits_{j\in N\setminus N^*}[\gamma_{ij}\cdot x_j + (1 - \gamma_{ij})\cdot a_{ij}]
\\&=&\sum \limits_{j\in N^*} a_{ij}+\sum \limits_{j\in N\setminus N^*, d_{p_jj} \leq a_{ij}}\gamma_{ij}\cdot x_j + \sum \limits_{j\in N\setminus N^*, a_{ij}\leq d_{(p_j-1)j}}(1 -\gamma_{ij})\cdot a_{ij}
\\&=&\sum \limits_{j\in N^*} a_{ij}+\sum \limits_{j\in N\setminus N^*, d_{p_jj} \leq a_{ij}} x_j+ \sum \limits_{j\in N\setminus N^*, a_{ij}\leq d_{(p_j-1)j}}a_{ij}
\\&=&\sum \limits_{j\in N^*} a_{ij}\wedge x_j+\sum \limits_{j\in N\setminus N^*, d_{p_jj} \leq a_{ij}}a_{ij}\wedge x_j+ \sum \limits_{j\in N\setminus N^*, a_{ij}\leq d_{(p_j-1)j}}a_{ij}\wedge x_j
\\&=&\sum \limits_{j\in N^*} a_{ij}\wedge x_j+\sum \limits_{j\in N\setminus N^*} a_{ij}\wedge x_j
\\&=&\sum \limits_{j\in N} a_{ij}\wedge x_j.
 \end{eqnarray*}
Since $b_i\leq \sum \limits_{j\in N^*} a_{ij}+\sum \limits_{j\in N\setminus N^*}[\gamma_{ij}\cdot x_j + (1 - \gamma_{ij})\cdot a_{ij}]=\sum \limits_{j\in N} a_{ij}\wedge x_j$ and $\lambda x_i=\sum \limits_{j\in N^*} a_{ij}+\sum \limits_{j\in N\setminus N^*}[\gamma_{ij}\cdot x_j + (1 - \gamma_{ij})\cdot a_{ij}]=\sum \limits_{j\in N} a_{ij}\wedge x_j$ for any $i\in N$, we have $x\in V^*(A, \lambda)$.
\end{proof}

Noting that in Theorem \ref{th4.1}, $\lambda$ can be seen as a given constant with $\lambda\in [0, +\infty)$. So that every solution $x= (x_1, x_2, \cdots, x_n)$ of system (\ref{ac}) can be represented by $\lambda$, and the corresponding $\lambda$ can be determined by the $x= (x_1, x_2, \cdots, x_n)$ since $x_j= d_{0j}=1\mbox{ for all } j\in N^*$ and $d_{{(p_j-1)}j} \leq x_j \leq d_{p_jj}\mbox{ for all } j\in N\setminus N^*$.

The proof of Theorem \ref{th4.1} implies the following two theorems.
\begin{theorem}\label{th4.3}
If system (\ref{ac}) is unsolvable for any $p=(p_1, p_2, \cdots, p_n)\in P$, then $V^*(A, \lambda)=\emptyset$.
\end{theorem}
\begin{theorem}\label{th4.2}
For any $p=(p_1, p_2, \cdots, p_n)\in P$, if $x= (x_1, x_2, \cdots, x_n)$ satisfies system (\ref{ac}), then $x\in V^*(A, \lambda)$.
\end{theorem}

For any $p=(p_1, p_2, \cdots, p_n)\in P$, denote the solution set of system (\ref{ac}) corresponding to $p=(p_1, p_2, \cdots, p_n)\in P$ by $V^*(A, \lambda, p)$. Then from Theorems \ref{th4.1} and \ref{th4.2}, the constrained eigenproblems of addition-min algebras are equivalent to solving system (\ref{ac}). So that we can summarize an algorithm to find all the constrained eigenvectors and eigenvalues of $A$ as follows.

\begin{algorithm}\label{al2}
Input $A=(a_{ij})_{n\times n}$ and $b$. Output $\Lambda^*(A)$ and $V^*(A)$.

Step 1. If $(1, 1, \cdots, 1)$ isn't a solution of system (\ref{b}) then $\Lambda^*(A)=\emptyset$, $V^*(A)=\emptyset$ and stop.

Step 2. Compute $P = \{(p_1, p_2, \cdots, p_n) | p_j\in P_j \mbox{ for any }j\in N\}$ defined by (\ref{p}).

Step 3. For any $p=(p_1, p_2, \cdots, p_n)\in P$, construct the corresponding system (\ref{ac}).

Step 4. Compute $V^*(A, \lambda, p)$ by the corresponding system (\ref{ac}).

Step 5. Output $\Lambda^*(A)=\bigcup \limits_{p\in P}\{\lambda\in [0,+\infty)|V^*(A,\lambda, p)\}$ and $V^*(A)=\bigcup\limits_{\lambda\in \Lambda^*(A)}V^*(A,\lambda)$.

Step 6. End.\end{algorithm}

\begin{theorem}\label{th4.3} Algorithm \ref{al2} terminates after $O((n+1)^{n})$ operations.
\end{theorem}
\begin{proof}Similar to the proof of Theorem \ref{th3.3}, one can prove that Algorithm \ref{al2} terminates after $O((n+1)^{n})$ operations since Step 3 is the key process of Algorithm \ref{al2} and $|P|\leq (n+1)^{n}$.
\end{proof}

 The following example illustrates Algorithm \ref{al2}.
\begin{example}\label{ex2}
\emph{Consider the following system:}
$$\left\{ \begin{array}{l}
0.4\wedge x_1 + 0.6\wedge x_2\geq 0.8,\\
0.2\wedge x_1 + 0.5\wedge x_2\geq 0.5,\\
0.4\wedge x_1 + 0.6\wedge x_2=\lambda x_1,\\
0.2\wedge x_1 + 0.5\wedge x_2=\lambda x_2.\\
\end{array} \right.$$
\end{example}
Step 1. $(1, 1, \cdots, 1)$ is a solution of system (\ref{b}).\\
Step 2. By $\check{\alpha}_1=0.2$ and $\check{\alpha}_2=0.4$. we have
$P_1=\{1, 2\}$ and $P_2=\{1, 2, 3\}$.
\begin{eqnarray*}
P=\{(1, 1), (1, 2),(1, 3), (2, 1), (2, 2), (2, 3)\}.
\end{eqnarray*}\\
Step 3. Construct and solve the corresponding system (\ref{ac}):
\begin{enumerate}
\item [(1)] For $p_1=(1, 1)$,
$$\left\{ \begin{array}{l}
x_1 + x_2 \geq 0.8,\\
0.2 + x_2\geq 0.5,\\
x_1 + x_2 =\lambda x_1,\\
0.2 + x_2=\lambda x_2,\\
0.2\leq x_1\leq 0.4,\\
0.4\leq x_2\leq 0.5.
\end{array} \right.$$
We get $x_1=\frac{0.2}{(\lambda-1)^2}$, $x_2=\frac{0.2}{\lambda-1}$ and $\lambda-1>0$. Since $0.2\leq \frac{0.2}{(\lambda-1)^2}\leq 0.4$, $0.4\leq \frac{0.2}{\lambda-1}\leq 0.5$ and $\frac{0.2}{(\lambda-1)^2}+\frac{0.2}{\lambda-1}\geq 0.8$, it is impossible. Therefore, $V^*(A, \lambda, p_1)=\emptyset.$

\item [(2)] For $p_2=(1, 2)$,
$$\left\{ \begin{array}{l}
x_1 + x_2 \geq 0.8,\\
0.2 + 0.5\geq 0.5,\\
x_1 + x_2 =\lambda x_1,\\
0.2 + 0.5=\lambda x_2,\\
0.2\leq x_1\leq 0.4,\\
0.5\leq x_2\leq 0.6.
\end{array} \right.$$
We get $x_1=\frac{0.7}{\lambda(\lambda-1)}$, $x_2=\frac{0.7}{\lambda}$ and $\lambda-1>0$. Since $0.2\leq \frac{0.7}{\lambda(\lambda-1)}\leq 0.4$, $0.5\leq \frac{0.7}{\lambda}\leq 0.6$ and $\frac{0.7}{\lambda(\lambda-1)}+\frac{0.7}{\lambda}\geq 0.8$, it is impossible. Therefore, $V^*(A, \lambda, p_2)=\emptyset.$

\item [(3)]  For $p_3=(1, 3)$,
$$\left\{ \begin{array}{l}
x_1 + 0.6 \geq 0.8,\\
0.2 + 0.5\geq 0.5,\\
x_1 + 0.6 =\lambda x_1,\\
0.2 + 0.5=\lambda x_2,\\
0.2\leq x_1\leq 0.4,\\
0.6\leq x_2\leq 1.
\end{array} \right.$$
We get $x_1=\frac{0.6}{\lambda-1}$, $x_2=\frac{0.7}{\lambda}$ and $ \lambda-1>0$. Since $0.2\leq \frac{0.6}{\lambda-1}\leq 0.4$ and $0.6\leq \frac{0.7}{\lambda}\leq 1$, it is impossible. Therefore, $V^*(A, \lambda, p_3)=\emptyset.$

\item [(4)] For $p_4=(2, 1)$,
$$\left\{ \begin{array}{l}
0.4 + x_2 \geq 0.8,\\
0.2 + x_2\geq 0.5,\\
0.4 + x_2 =\lambda x_1,\\
0.2 + x_2=\lambda x_2,\\
0.4\leq x_1\leq 1,\\
0.4\leq x_2\leq 0.5.
\end{array} \right.$$
We get $x_1=\frac{0.4\lambda-0.2}{\lambda(\lambda-1)}$, $x_2=\frac{0.2}{\lambda-1}$ and $\lambda-1>0$. Since $0.4\leq \frac{0.4\lambda-0.2}{\lambda(\lambda-1)}\leq 1$ and $0.4\leq \frac{0.2}{\lambda-1}\leq 0.5$, we have $\frac{7}{5}\leq\lambda\leq \frac{3}{2}$. Therefore, $V^*(A, \lambda, p_4)=\{(\frac{0.4\lambda-0.2}{\lambda(\lambda-1)}, \frac{0.2}{\lambda-1})|\frac{7}{5}\leq\lambda\leq \frac{3}{2}\}.$

\item [(5)] For $p_5=(2, 2)$,
$$\left\{ \begin{array}{l}
0.4 + x_2 \geq 0.8,\\
0.2 + 0.5\geq 0.5,\\
0.4 + x_2 =\lambda x_1,\\
0.2 + 0.5=\lambda x_2,\\
0.4\leq x_1\leq 1,\\
0.5\leq x_2\leq 0.6.
\end{array} \right.$$
We get $x_1=\frac{0.4\lambda+0.7}{\lambda^2}$, $x_2=\frac{0.7}{\lambda}$ and $\lambda>0$. Since $0.4\leq \frac{0.4\lambda+0.7}{\lambda^2}\leq 1$ and $0.5\leq \frac{0.7}{\lambda}\leq 0.6$, we have $\frac{7}{6}\leq\lambda\leq \frac{7}{5}$. Therefore, $V^*(A, \lambda, p_5)=\{(\frac{0.4\lambda+0.7}{\lambda^2}, \frac{0.7}{\lambda})|\frac{7}{6}\leq\lambda\leq \frac{7}{5})\}.$

\item [(6)] For $p_6=(2, 3)$,
$$\left\{ \begin{array}{l}
0.4 + 0.6\geq 0.8,\\
0.2 + 0.5\geq 0.5,\\
0.4 + 0.6 =\lambda x_1,\\
0.2 + 0.5=\lambda x_2,\\
0.4\leq x_1\leq 1,\\
0.6\leq x_2\leq 1.
\end{array} \right.$$
We get $x_1=\frac{1}{\lambda}$, $x_2=\frac{0.7}{\lambda}$ and $\lambda>0$. Since $0.4\leq \frac{1}{\lambda}\leq 1$ and $0.6\leq \frac{0.7}{\lambda}\leq 1$, we have $1\leq\lambda\leq \frac{7}{6}$. Therefore, $V^*(A, \lambda, p_6)=\{(\frac{1}{\lambda}, \frac{0.7}{\lambda})|1\leq\lambda\leq \frac{7}{6}\}.$

\end{enumerate}
Step 4. Output
\begin{eqnarray*}
 \Lambda^*(A)&=&\bigcup \limits_{p\in P}\{\lambda\in [0,+\infty)|V^*(A,\lambda, p)\}
 \\&=&\{\lambda\in [0,+\infty)|V^*(A,\lambda, p_4)\}\cup \{\lambda\in [0,+\infty)|V^*(A,\lambda, p_5)\}\\&&\cup\{\lambda\in [0,+\infty)|V^*(A,\lambda, p_6)\}
 \\&=&\{\lambda|\frac{7}{5}\leq\lambda\leq \frac{3}{2}\}\cup\{\lambda|\frac{7}{6}\leq\lambda\leq \frac{7}{5})\}\cup\{\lambda|1\leq\lambda\leq \frac{7}{6}\}
 \\&=&\{\lambda|1\leq\lambda\leq\frac{3}{2}\}
\end{eqnarray*}
and
\begin{eqnarray*}
 V^*(A)&=&\bigcup \limits_{\lambda\in\Lambda^*(A)}V^*(A, \lambda)
 \\&=&\{(\frac{0.4\lambda-0.2}{\lambda(\lambda-1)}, \frac{0.2}{\lambda-1})|\frac{7}{5}\leq\lambda\leq \frac{3}{2}\}\cup \{(\frac{0.4\lambda+0.7}{\lambda^2}, \frac{0.7}{\lambda})|\frac{7}{6}\leq\lambda\leq \frac{7}{5})\}\\&&\cup\{(\frac{1}{\lambda}, \frac{0.7}{\lambda})|1\leq\lambda\leq \frac{7}{6}\}.
\end{eqnarray*}

\section{Supereigenproblems of addition-min algebras}\label{se5}

In order to improve the enthusiasm of the terminals, Yang et al. \cite{Yang2019} considered the ratio of the data-download quality to the data-sending quality, i.e., they introduced a supereigenproblem of addition-min algebra
\begin{equation}\label{equ01}A \odot x^T \geq \lambda x^T
\end{equation}
where $A=(a_{ij})_{n\times n}$ with $a_{ii}=0$ for all $i\in N$. They obtained a finite number of the supereigenvectors of system \eqref{equ01} associated with a given supereigenvalue $\lambda\in (0, n-1]$. They further investigated a constrained supereigenproblem of addition-min algebra, more specific speaking, they studied the maximum constrained supereigenvalue and the corresponding constrained supereigenvector under system (\ref{b}), i.e.,
 \begin{align}\nonumber
&\mbox{max}\quad \lambda\\
&\mbox{s.t.}\quad\left\{ \begin{array}{l}\label{equ02}
A \odot x^T \geq b^T,\\
A \odot x^T \geq \lambda x^T
\end{array} \right.
\end{align}
where $A=(a_{ij})_{n\times n}$ with $a_{ii}=0$ for all $i\in N$. They developed a nonlinear programming approach to find the unique maximum constrained supereigenvalue and one corresponding constrained supereigenvector of system (\ref{equ02}). In this section, we make Algorithm \ref{al1} be suitable for characterizing the feasible region of all the supereigenvectors of system \eqref{equ01} associated with a given supereigenvalue $\lambda\in (0, n-1]$, and present an algorithm to obtain the maximum constrained supereigenvalue and the feasible region of all the corresponding constrained supereigenvectors of system (\ref{equ02}).

Just replacing $A \odot x^T = \lambda x^T$ by $A \odot x^T \geq \lambda x^T$, and $'='$ in the corresponding content of Section \ref{se3} by $'\geq'$, one can easily prove that both the corresponding Theorems \ref{th3.1} and \ref{th3.2} hold, and Algorithm \ref{al1} is suitable for describing the feasible region of all the supereigenvectors of system \eqref{equ01} associated with a given supereigenvalue $\lambda\in (0, n-1]$. We omit the detail and only use the following example to illustrated it.
 \begin{example}\label{ex3}
\emph{Consider the feasible region of all the supereigenvectors of the following system with $\lambda=1$:}
$$\left\{ \begin{array}{l}
0\wedge x_1 + 0.6\wedge x_2\geq x_1,\\
0.4\wedge x_1 + 0\wedge x_2\geq x_2.
\end{array} \right.$$
\end{example}
Step 1. Compute $K_1=\{1, 2\}$ and $K_2=\{1, 2\}$. Then
\begin{eqnarray*}
K=\{(1, 1), (1, 2), (2, 1), (2, 2)\}.
\end{eqnarray*}
Step 2. Construct and solve the corresponding system (\ref{ab}) (replacing $'='$ of system (\ref{ab}) by $'\geq'$):
\begin{enumerate}
\item [(1)] For $k_1=(1, 1)$,
$$\left\{ \begin{array}{l}
x_2 \geq x_1,\\
x_1 \geq x_2,\\
0\leq x_1\leq 0.4,\\
0\leq x_2\leq 0.6.
\end{array} \right.$$
We get $x_1=x_2$. Let $x_1=x_2=t$. Since $0 \leq x_1\leq 0.4$ and $0 \leq x_2\leq 0.6$, we have $t\in (0, 0.4]$. Therefore, the feasible region in this case is $\{(t,t)|t\in (0, 0.4]\}$.

\item [(2)] For $k_2=(1, 2)$,
$$\left\{ \begin{array}{l}
0.6 \geq x_1,\\
x_1\geq x_2,\\
0\leq x_1\leq 0.4,\\
0.6\leq x_2\leq 1.
\end{array} \right.$$
Obviously, it is impossible.

\item [(3)] For $k_3=(2, 1)$,
$$\left\{ \begin{array}{l}
x_2\geq x_1,\\
0.4 \geq x_2,\\
0.4\leq x_1\leq 1,\\
0\leq x_2\leq 0.6.
\end{array} \right.$$
We get $x_1=0.4, x_2=0.4$. Then the feasible region in this case is $\{(0.4, 0.4)\}$.

\item [(4)] For $k_4=(2, 2)$,
$$\left\{ \begin{array}{l}
0.6 \geq x_1,\\
0.4\geq x_2,\\
0.4\leq x_1\leq 1,\\
0.6\leq x_2\leq 1.
\end{array} \right.$$
Obviously, it is impossible.
\end{enumerate}
Step 3. The feasible region of all the supereigenvectors associated with $\lambda=1$ is $\{(t,t)|t\in (0, 0.4]\}$.

In what follows, we suggest an algorithm to obtain the maximum constrained supereigenvalue and the feasible region of all the corresponding constrained supereigenvectors of system (\ref{equ02}).

 From the proof of Theorem \ref{th4.1}, for any solution $x$ of system \eqref{equ02} there exists a corresponding $p=(p_1, p_2, \cdots, p_n)\in P$ such that
 \begin{align*}
\left\{ \begin{array}{l}
x_j= d_{0j}=1\mbox{ for all } j\in N^*,\\
d_{{(p_j-1)}j} \leq x_j \leq d_{p_jj}\mbox{ for all } j\in N\setminus N^*,\\
\sum \limits_{j\in N} a_{ij}\wedge x_j=\sum \limits_{j\in N^*}a_{ij}+\sum \limits_{j\in N\setminus N^*}[\gamma_{ij}\cdot x_j + (1 - \gamma_{ij})\cdot a_{ij}]\mbox{ for all } i\in N.
\end{array} \right.
\end{align*}
 Thus system \eqref{equ02} corresponding to $p=(p_1, p_2, \cdots, p_n)\in P$ can be written as follows.
\begin{align}\nonumber
&\max\quad\lambda\\
&\mbox{s.t.}\quad\left\{ \begin{array}{l}\label{ae}
x_j= d_{0j}=1\mbox{ for all } j\in N^*,\\\\
d_{{(p_j-1)}j} \leq x_j \leq d_{p_jj}\mbox{ for all } j\in N\setminus N^*,\\ \\
\sum \limits_{j\in N^*}a_{ij}+\sum \limits_{j\in N\setminus N^*}[\gamma_{ij}\cdot x_j + (1 - \gamma_{ij})\cdot a_{ij}]\geq b_i\mbox{ for all } i\in N,\\\\
\sum \limits_{j\in N^*}a_{ij}+\sum \limits_{j\in N\setminus N^*}[\gamma_{ij}\cdot x_j + (1 - \gamma_{ij})\cdot a_{ij}]\geq\lambda x_i\mbox{ for all } i\in N
\end{array} \right.
\end{align}
where the operation $``\cdot"$ represents the ordinary multiplication and
$$\gamma_{ij} = \left\{ \begin{array}{l}
1, d_{p_jj} \leq a_{ij},\\
0, a_{ij} \leq d_{(p_j-1)j}.
\end{array} \right.$$

For any $p=(p_1, p_2, \cdots, p_n)\in P$, denote the local optimal maximum constrained supereigenvalue of system \eqref{ae} by $\lambda(p)$ and the corresponding feasible region of the constrained supereigenvectors by $V^0(A, \lambda(p), p)$. Then the maximum constrained supereigenvalue $\lambda=\max\{\lambda(p)|p\in P\}$ and the feasible region of all the corresponding constrained supereigenvectors $V^0(A, \lambda)=\bigcup\limits_{\lambda=\lambda(p)}V^0(A, \lambda(p), p)$ where $\lambda(p)$ can be computed by using the software LINGO or MATLAB. We can summarize an algorithm to find the maximum constrained supereigenvalue and the feasible region of all the corresponding constrained supereigenvectors of $A$ as follows.
\begin{algorithm}\label{al3}
Input $A=(a_{ij})_{n\times n}$ and $b$. Output $\lambda$ and $V^0(A, \lambda)$.

Step 1.  If $(1, 1, \cdots, 1)$ isn't a solution of system (\ref{b}) then $\lambda$ doesn't exist, $V^0(A, \lambda)=\emptyset$ and stop.

Step 2. Compute $P = \{(p_1, p_2, \cdots, p_n) | p_j\in P_j \mbox{ for any }j\in N\}$ defined by (\ref{p}).

Step 3. For any $p=(p_1, p_2, \cdots, p_n)\in P$, construct the corresponding system (\ref{ae}).

Step 4. Solve the corresponding system (\ref{ae}), and obtain $\lambda(p)$ by the software LINGO or MATLAB and $V^0(A, \lambda(p), p)$.

Step 5. Output $\lambda=\max\{\lambda(p)|p\in P\}$ and $V^0(A, \lambda)=\bigcup\limits_{\lambda=\lambda(p)}V^0(A, \lambda(p), p)$.

Step 6. End.
\end{algorithm}

Notice that similar to Theorem \ref{th4.3}, one can see that Algorithm \ref{al3} terminates after $O((n+1)^{n})$ operations.

\begin{example}\label{ex4}
\emph{Consider the following problem:}
\begin{align*}
&\emph{max}\quad \lambda\\
&\emph{s.t.}\quad\left\{ \begin{array}{l}
0\wedge x_1 + 0.6\wedge x_2\geq 0.2,\\
0.4\wedge x_1 + 0\wedge x_2\geq 0.3,\\
0\wedge x_1 + 0.6\wedge x_2\geq\lambda x_1,\\
0.4\wedge x_1 + 0\wedge x_2\geq\lambda x_2.
\end{array} \right.
\end{align*}
\end{example}
Step 1. $(1, 1, \cdots, 1)$ is a solution of system (\ref{b}).\\
Step 2. By $\check{\alpha}_1=0$ and $\check{\alpha}_2=0$, we have $P_1=\{1, 2\}$ and $P_2=\{1, 2\}$. Then
\begin{eqnarray*}
P=\{(1, 1), (1, 2), (2, 1), (2, 2)\}.
\end{eqnarray*}
Step 3. Construct and solve the corresponding system (\ref{ae}):
\begin{enumerate}
\item [(1)] For $p_1=(1, 1)$,
\begin{align*}
&\max\quad \lambda\\
&\mbox{s.t.}\quad\left\{ \begin{array}{l}
x_2 \geq 0.2,\\
x_1 \geq 0.3,\\
x_2 \geq \lambda x_1,\\
x_1 \geq \lambda x_2,\\
0\leq x_1\leq 0.4,\\
0\leq x_2\leq 0.6.
\end{array} \right.
\end{align*}
By MATLAB, $\lambda(p_1)=1$. Then $x_1=x_2$. Let $x_1=x_2=t$. Since $0.3\leq x_1\leq 0.4$ and $0.2\leq x_2\leq 0.6$, we have $0.3\leq t\leq 0.4$. Therefore, $V^0(A, 1, p_1)=\{(t, t)|t\in [0.3, 0.4]\}$.

\item [(2)] For $p_2=(1, 2)$,
\begin{align*}
&\mbox{max}\quad \lambda\\
&\mbox{s.t.}\quad\left\{ \begin{array}{l}
0.6 \geq 0.2,\\
x_1\geq 0.3,\\
0.6 \geq \lambda x_1,\\
x_1\geq \lambda x_2,\\
0\leq x_1\leq 0.4,\\
0.6\leq x_2\leq 1.
\end{array} \right.
\end{align*}
By MATLAB, $\lambda(p_2)=\frac{2}{3}$. Then $0.6 \geq \frac{2}{3}x_1$ and $x_1\geq \frac{2}{3}x_2$. Since $0.3\leq x_1\leq 0.4$ and $0.6\leq x_2\leq 1$, we have $x_1=0.4$ and $x_2=0.6$. Therefore, $V^0(A, \frac{2}{3}, p_2)=\{(0.4, 0.6)\}$.

\item [(3)] For $p_3=(2, 1)$,
\begin{align*}
&\mbox{max}\quad \lambda\\
&\mbox{s.t.}\quad\left\{ \begin{array}{l}
x_2\geq 0.2,\\
0.4 \geq 0.3,\\
x_2\geq \lambda x_1,\\
0.4 \geq \lambda x_2,\\
0.4\leq x_1\leq 1,\\
0\leq x_2\leq 0.6.
\end{array} \right.
\end{align*}
By MATLAB, $\lambda( p_3)=1$. Then $x_2\geq x_1$ and $0.4 \geq x_2$. Since $0.4\leq x_1\leq 1$ and $0\leq x_2\leq 0.4$, we have $x_1=x_2=0.4$. Therefore, $V^0(A, 1, p_3)=\{(0.4,0.4)\}$.

\item [(4)] For $p_4=(2, 2)$,
\begin{align*}
&\mbox{max}\quad \lambda\\
&\mbox{s.t.}\quad\left\{ \begin{array}{l}
0.6 \geq 0.2,\\
0.4\geq 0.3,\\
0.6 \geq \lambda x_1,\\
0.4\geq \lambda x_2,\\
0.4\leq x_1\leq 1,\\
0.6\leq x_2\leq 1.
\end{array} \right.
\end{align*}
By MATLAB, $\lambda(p_4)=\frac{2}{3}$. Then $0.6 \geq \frac{2}{3}x_1$ and $0.4\geq \frac{2}{3}x_2$. Since $0.4\leq x_1\leq 1$ and $0.6\leq x_2\leq 1$, we have $0.4\leq x_1\leq0.9$ and $x_2=0.6$. Let $x_1=t$. Then $V^0(A, \frac{2}{3}, p_4)=\{(t,0.6)|t\in [0.4, 0.9]\}$.
\end{enumerate}
Step 4. Output $$\lambda=\max\{\lambda(p)|p\in P\}=1$$ and $$V^0(A, 1)=\bigcup\limits_{\lambda(p)=1}V^0(A, \lambda(p), p)=V^0(A, 1, p_1)\cup V^0(A, 1, p_3)= \{(t,t)|t\in [0.3, 0.4]\}.$$
\section{Concluding remark}\label{se6}
In this article, we first investigated the eigenproblems of addition-min algebras and showed Algorithm \ref{al1} for finding all the eigenvalues and eigenvectors for a given matrix. Then, we studied the constrained eigenproblems of addition-min algebras and proposed Algorithm \ref{al2} for computing all the constrained eigenvectors and eigenvalues for a given matrix. We finally discussed the supereigenproblems of addition-min algebras and suggest Algorithm \ref{al3} for obtaining the maximum constrained supereigenvalue and depicting the feasible region of all the constrained supereigenvectors for a given matrix. Since the computational complexity of our algorithms is $O((n+1)^n)$, they may involve heavy and complicated work when $n$ is a larger number. Therefore, finding an efficient method for the eigenproblems (resp. the constrained eigenproblems and the supereigenproblems) of addition-min algebras is an interesting problem in the future.

\end{document}